\documentclass[12pt]{article}
\usepackage{amssymb,amsfonts,euscript}
\usepackage{amsmath}

\newcounter{minutes}
\divide\time by 60
\newcounter{hours}
\multiply\time by 60
\addtocounter{minutes}{-\time}

\begin{document}
\begin{center}
 \LARGE\bf  Robin functions and distortion theorems\\ 
for regular mappings\\
\end{center}
\begin{center}
\texttt{File:~\jobname .tex, 2007-01-23, 
        printed: \number\year-\number\month-\number\day, 
        \thehours.\ifnum\theminutes<10{0}\fi\theminutes}
\end{center}

\begin{center}\large\bf  V.N. Dubinin and M. Vuorinen
\end{center}

\vspace{0.7cm}  {\bf Abstract.} Capacities of generalized 
condensers are applied to prove a two-point distortion theorem 
for conformal mappings. The result is expressed in terms of the 
Robin function and the Robin capacity with respect to the domain of 
definition of the mapping and subsets of the boundary of this domain. 
The behavior of Robin function 
under multivalent functions is studied. Some corollaries and examples of 
applications to distortion theorems for regular functions are given.

\vspace{0.7cm} \noindent {\bf  1 Introduction}

\vspace{0.5cm}\noindent The notion of the Robin function arises in a natural 
way in the theory of partial differential equations as a generalization 
of Green's function [BS]. Because of conformal invariance these functions have 
wide applicability in complex analysis. The study of Robin's function and 
the Robin capacity associated to it has recently received a lot of attention 
(see e.g. ([DS1], [DS2], [DPT], [D1], [DitSol], [NT], 
[Nas],[St1], [St2], [Vas]). 
However, a series of questions remains open. Thus, the behavior of the Robin 
capacity under multivalent mappings has not been adequately studied. 
Furthermore in the literature not enough attention has been paid to 
applications of the Robin capacity to distortion theorems of regular 
functions (even for the case of conformal mappings). The capacity approach 
in the proof of such results is not restricted merely to the case of plane 
but it also applies to the case of higher dimensions [V]. 
In this paper these gaps will be partially filled. We begin with applications 
to conformal maps. In Section 2 a new two-point distortion theorem for regular 
univalent functions is proved in terms of the Robin capacity and the 
Robin function (Theorem 2.1). 
This theorem has a general character and contains many earlier results as particular cases
(see Section 3). The proof of Theorem 2.1 makes use of the connection of some quadratic 
forms which depend on the values of the Robin function at the points under consideration 
and asymptotic expansion of the capacity of generalized condensers [Dub1]. 
Note that a considerable part of well-known two-point distortion 
theorems were obtained by methods applicable to simply-connected 
domains for which the study of functions 
defined for instance in a ring domain is already quite complicated. 
It should be emphasized that our method 
applies in the same way for an arbitrary number of points and for domains 
of arbitrary connectivity. 
In Section 4 the behavior of Robin functions under multivalent 
functions is studied. For this purpose Lindel\"of's principle 
involving the Green function 
is generalized to the case of the Robin function (Theorem 4.1). 
This generalization leads to Theorem 4.2 about 
the behaviour of the Robin capacity under regular mappings. Here we develop further the 
approach of Mityuk [M1,M2], from which many applications follow. For these applications
a crucial feature is taking into account 
the multiplicity of the covering. We also consider particular cases of 
these theorems and other results, connected with aforementioned questions. 
We proceed to the definition of the Robin function and to its counterpart 
with a boundary pole [Dub1].

Let a domain $B$ of the complex plane $\overline{\mathbf{C}}_z$ be 
bounded by a finite number of analytic curves, let $\gamma$ be a 
non-empty closed subset of $\partial B,$ consisting of a finite 
number of nondegenerated Jordan arcs, and let $z_0$ be a finite 
point of the set $\overline{B}\,\backslash\gamma.$ We denote by 
$g(z)=g_B(z, z_0, \gamma)$ a continuous real-valued function on 
$\overline{B}\,\backslash \{z_0\},$ continuously differentiable on 
$\overline{B}\,\backslash (\gamma\cup\{z_0\}),$ harmonic in 
$B\,\backslash \{z_0\}$ and satisfying the following conditions

 \vskip0.2cm
$\qquad\quad g(z)=0\,\,\,\,\,\,\qquad \qquad\qquad
\,\,\,\,\mbox{for} \qquad z \in \gamma\,\,\,, $

 \vskip0.4cm
$\qquad\quad{\displaystyle\frac{\partial g}{\partial
n}}(z)=0\qquad \quad\quad\, \quad\quad\,\,\,\mbox{for}
\qquad\, z \in (\partial B)\big\backslash\left(\gamma\cup
\{z_0\}\right)\,\,\,, $\\
$ \quad g(z)+ \log|z-z_0|$ is a harmonic function   
in a neighborhood of the point $z_0$ ($\partial\big/\partial n$ 
denotes differentiation in the direction of the interior normal 
to the boundary $\partial B$). In the case when $z_0= \infty,$ 
the function $g_B(z, z_0, \gamma)$ is defined in an analogical way 
with the only difference, that the harmonicity of the function 
$g_B(z, z_0, \gamma)-\log|z|$ in a neighborhood of infinity is required. 
A finitely connected domain $B \subset \overline{\mathbf{C}}_z$ and 
a closed subset $\gamma$ of its boundary $\partial B$ are called admissible, 
if the boundary $\partial B$ does not have isolated points, 
and $\gamma$ consists of a finite number of nondegerated boundary components. 
The definition of the function $g_B(z, z_0, \gamma)$ for general 
admissible domains $B$ and sets $\gamma$ takes place with the help of 
a conformal mapping (cf. [Dub1]). For $z_0 \in B$ the function 
$g_B(z, z_0, \gamma)$ is called the Robin function of the domain 
$B$ with a pole at the point $z_0.$ In the case when $\gamma=\partial B$, 
the Robin function agrees with the Green function 
$g_B(z, z_0, \partial B)= g_B(z, z_0).$ We also introduce the notation:
   \[
   r(B, \gamma, z_0)=\exp\{\lim\limits_{z\to z_0}[g_B(z, z_0,\gamma)+\log|z-z_0|]\}
   \]
   in the case of a finite point $z_0$ and
   \[
r(B, \gamma, \infty)=\exp\{\lim\limits_{z\to \infty}[g_B(z,
\infty,\gamma)+\log|z|]\}.
   \]
   When $z_0 \in B$ and $\gamma=\partial B$ the quantity $r(B, \gamma, z_0)\equiv r(B, z_0)$ 
is called the inner radius of the domain $B$  with respect to the point $z_0$. 
If furthermore $B$ is simply connected, then $r(B, z_0)$ is its conformal radius. 
For $z_0=\infty$ the conformal radius is sometimes defined as $r^{-1}(B, z_0).$ 
The quantity $r^{-1}(B,\infty)$ agrees with the logarithmic capacity of the complement 
$r^{-1}(B, \infty)=\mathrm{cap}\,(\overline{\mathbf{C}}_z\backslash B).$ 
For $z_0=\infty$ and $z_0 \in B$ the quantity $r^{-1}(B,\gamma, z_0)$ is called 
{\it the Robin capacity} of the set $\gamma$ with repect to the domain $B$.

Finally, when $B$ is simply connected, $\gamma$ is a boundary arc of 
$B$  and $z_0 \in \partial B\backslash\gamma ,$  the quantity 
$r(B, \gamma, z_0)$ has been considered under different names 
and from various viewpoints in   [Her], [GH], [Mik], [Sol1].

The main method for the proof of Theorem 2.1 is the notion of 
the capacity of a generalized condenser [Dub1] 
(in what follows the word "generalization" will be omitted). Let $B$
be a finitely connected domain on the complex plane 
$\overline{\mathbf{C}}_z$, and let $\overline{B}$ denote 
the compactification of $B$ by means of prime ends in the sense of 
Caratheodory. In what follows, whenever it makes sense we will 
identify an element of $\overline{B}$, corresponding to an interior point 
of $B$, with the same point, and the support of an accessible boundary point 
and this point itself will be denoted with the same letter. 
A triple $C=(B, \mathcal{E}, \Delta)$ is called a condenser where 
$\mathcal{E}=\{E_k\}^n_{k=1}$ is the union of pairwise nonintersecting 
subsets $E_k, \, k=1, \ldots, n$, closed in $\overline{B}$ and 
$\Delta=\{t_k\}^n_{k=1}$ is the union of real numbers 
$t_k,\,  k=1, \ldots, n$. The capacity $\mathrm{cap}\, C$ of 
the condenser $C$ is defined as the infimum of Dirichlet integrals
  \[
  I(v, B):=\int\int_B|\nabla v|^2dxdy
  \]
taken over all functions $v(z)\, (z=x+iy),$ continuous in 
$\overline{B}$, satisfying Lipschitz condition at some neighborhood 
of each finite point $B$ and equal to $t_k$ at some neighborhood of 
$E_k, \, k=1, \ldots, n$.

We also need the asymptotics of the capacity of a condenser, a 
special case of the result in [Dub1]. For a finite point $z_0$ of the 
complex sphere $\overline{\mathbf{C}}_z$ we denote by $E(z_0, r)$ a closed disk 
with the center at the point $z_0$ and with radius $r>0.$ In the case of the point at 
infinity we define $E(\infty, r):=\{z:|z|\geq 1/r\}.$ Let the domain 
$B\subset \overline{\mathbf{C}}_z$ and let $\gamma\subset\partial B$ be admissible, 
let $Z=\{z_k\}^n_{k=1}$ be the union of distinct points of the domain $B$, let 
$\triangle=\{t_k\}^n_{k=1}$ be the union of real numbers, 
$\sum\limits_{k=1}^n\, t_k^2\neq 0,$ and let $\Psi=\{\psi_k(r)\}^n_{k=1}$, where 
$\psi_k(r)\equiv\mu_k r^{\nu_k},$ and $  \mu_k, \,\nu_k, k=1, \ldots, n$ 
are positive numbers. For sufficiently small $r$ we denote\quad
$C(r;B, \gamma, Z,$ $\Delta, \Psi)$
$:=(B, \{\gamma,$
$E(z_1, \psi_1(r)),$
$ \ldots,$
$E(z_n, \psi_n(r))\},$
$\{0, t_1,$
$\ldots, t_n\})$.
From Theorem 7 of [Dub1] the following asymptotic formula is obtained:

  \[
\mathrm{cap}\, C(r;B, \gamma, Z, \Delta,
\Psi)=2\pi\left[\sum\limits_{k=1}^n\,
\frac{t^2_k}{\nu_k}\right]\left(-\frac{1}{\log r}\right)-
\]

\[
-2\pi\left[\sum\limits_{k=1}^n\,
\frac{t^2_k}{\nu_k^2}\log\frac{r(B, \gamma,
z_k)}{\mu_k}+\sum\limits_{k=1}^n\sum\limits_{l=1\atop l\neq k}^n\,
\frac{t_k}{\nu_k} \frac{t_k}{\nu_l} g_B(z_k, z_l,
\gamma)\right]\left(\frac{1}{\log r}\right)^2+
\]

\[
\qquad\qquad\qquad\qquad \qquad +o\left(\left(\frac{1}{\log
r}\right)^2\right), \,\,\,\, r \rightarrow 0. \qquad
\qquad\qquad\qquad(1.1)
\]

Formula (1.1) remains valid if in the definition of the condenser 
$C(r;$ $B,$ $\gamma,$ $Z,$ $\Delta, \Psi)$ the disks 
$E(z_k, \psi_k(r)), k=1, \ldots, n$ are replaced with "almost disks" [Dub1].


 \vspace{2.0cm}\noindent {\Large\bf 2 Two-point distortion theorem for univalent functions}

\vspace{0.5cm}\noindent 
As well-known, 
quite many recent papers deal with two-point distortion theorems 
(cf. for example [J],[KMin],[MaMin],[KrR] and the literature therein). 
Our approach is different because we consider multiply connected domains 
and make use of the behavior of the function on the boundary. 
We first define "the weighted derivative" of a conformal mapping. 
Let the domain $B$ of the plane $\overline{\mathbf{C}}_z$ and the set 
$\gamma\subset\partial B$ be admissible. In analogy with 
the hyperbolic case we consider the conformal invariant
\[
\delta(\zeta, z; B,\gamma):=\exp\{-g_B(\zeta, z, B,\gamma)\},
\]
where the points $\zeta, z\in B$. Let the function $w=f(z)$ map 
conformally and univalently the domain $B$ into an admissible 
domain $G \subset \overline{\mathbf{C}}_w$ and let 
$\Gamma\subset\partial G$ be also admissible. For each point $z \in B$ we set
\[
|Df(z)|:=\lim\limits_{\zeta\to z}\,\,\frac{\delta(f(\zeta), f(z);
G, \Gamma)}{\delta(\zeta, z; B, \gamma)}=\frac{r(B, \gamma,
z)}{r(G, \Gamma, f(z))}\, |f'(z)|.
\]
This definition extends to the case of boundary points 
$z \in (\partial B)\backslash\gamma$ $(f(z) \in(\partial G)\backslash\Gamma)$, 
if the derivative $f'(z)$ is understood 
for example,
as the angular derivative. In the particular case where 
$B=\{z:|z|<1\},$ $\gamma=\partial B$ and $G=\{w:|w|<1\},$
$\Gamma=\partial G$ we have
\[
\delta(\zeta, z; B, \gamma
)=\left|\frac{\zeta-z}{1-\overline{z}\zeta}\right|
\]
with a similar formula for $\delta(f(\zeta), f(z); G, \Gamma)$ and,
therefore, the expression $|Df(z)|$ agrees with 
"hyperbolic--hyperbolic derivative"\, of $f$ at $z \in B$
\[
|Df(z)|=\frac{1-|z|^2}{1-|f(z)|^2}\,\,|f'(z)|
\]
[KrR, p. 116]. The following two-point distortion theorem contains, 
as special cases, several well-known results of this kind.

\vspace{0.5cm} {\bf  2.1. Theorem.}  {\it Let domains} $B$ {\it and} $G$, 
{\it which are subsets of } $\overline{\mathbf{C}}_z$ {\it and} 
$\overline{\mathbf{C}}_w$, {\it respectively, and the sets} $
\gamma\subset\partial B$ {\it and} $\Gamma\subset\partial G$ 
{\it be admissible. We assume, that the function} $f$ 
{\it maps the domain} $B$ {\it conformally and univalently into the domain} 
$G$ {\it such that} $f(\gamma)\subset\Gamma$, {\it where the image} 
$f(\gamma)$ {\it is understood in the sense of boundary correspodence. 
Then for all points} $z_1, z_2 \in B$ {\it and all real numbers} 
$t_1$  {\it and} $t_2$ {\it the following inequality holds}
\[
 \qquad\quad|D f(z_1)|^{t_1^2}
|D f(z_2)|^{t_2^2} \geq \left[\frac{\delta(z_1, z_2; B,
\gamma)}{\delta(f(z_1), f(z_2); G, \Gamma)
}\right]^{2t_1t_2}\,\,\, . \quad\qquad\quad (2.1)
\]
{\it If for a conformal mapping} $f$ {\it the inclusion} 
$f(B) \subset G$ {\it is valid and} $f((\partial B)\backslash\gamma)\subset
(\partial G)\backslash\Gamma$, {\it then}
\[
 \qquad|D f(z_1)|^{t_1^2}
|D f(z_2)|^{t_2^2} \leq \left[\frac{\delta(z_1, z_2; B,
\gamma)}{\delta(f(z_1), f(z_2); G, \Gamma) }\right]^{2t_1t_2}
\quad\qquad\quad\quad (2.2)
\]
{\it for all points} $z_1, z_2 \in B$ {\it and all real numbers}
$t_1, t_2.$

\vspace{0.5cm} {\bf  Proof.} We may assume that $t_1^2+ t_2^2\neq 0,$ 
and the domains $B$ and $G$ are bounded by analytic Jordan arcs. 
Suppose that the inclusion $f(\gamma)\subset\Gamma$ holds. 
For sufficiently small $r >0$ we consider the condenser
\[
C(r;f(B), f(\gamma), W, \Delta, \Psi),
\]
where $W=\{w_k\}^2_{k=1}$, $w_k=f(z_k), \, k=1, 2;$ 
$\Delta=\{t_1, t_2\};$ $\Psi=\{|f'(z_k)|r\}^2_{k=1}$. 
We shall prove that the capacity of this condenser does not 
exceed the capacity of the condenser
\[
C(r; G, \Gamma,  W, \Delta, \Psi).
\]
Indeed, if a function $v$ is continuous in $\overline{{G}}$, 
satisfies Lipschitz condition in a neighborhood of every finite 
point of the domain $G$, is equal to zero in a neighborhood 
of $\Gamma$ and equal to $t_k$ in a neighborhood of the set 
$E(w_k, |f'(z_k)|r),\, k=1, 2$, then it enjoys the aforelisted properties, 
if the domain $G$ is replaced with $f(B)$, and the set $\Gamma$ with $f(\gamma).$ Thus
\[
I(v, G)\geq I(v,f(B))\geq\mathrm{cap}\, C(r;f(B), f(\gamma), W,
\Delta, \Psi).
\]
Taking here the infimum over all possible functions $v$, we obtain
\[
\mathrm{cap}\, C(r;G, \Gamma, W, \Delta, \Psi)\geq\mathrm{cap}\,
C(r;f(B), f(\gamma), W, \Delta, \Psi).
\]
On the other hand, from the conformal invariance of the capacity it follows that
\[
\mathrm{cap}\, C(r;f(B), f(\gamma), W, \Delta, \Psi)=
\mathrm{cap}(B, \mathcal{E}_1, \Delta_1),
\]
where 
\[ \mathcal{E}_1=\{\gamma, f^{-1}(E(w_1, |f'(z_1)|r)),
f^{-1}(E(w_2, |f'(z_2)|r))\},\]
$\Delta_1=\{0, t_1, t_2\}.$ The sets $f^{-1}(E(w_k, |f'(z_k)|r)),$  $k=1,2$, 
form "almost disks". By virtue of formula (1.1)
\[
\mathrm{cap}\,(B, \mathcal{E}_1,
\Delta_1)=2\pi(t_1^2+t_2^2)\left(-\frac{1}{\log
r}\right)-2\pi\left[t_1^2\log r(B, \gamma, z_1)+\right.
\]
\[
\left.+t_2^2\log r(B, \gamma, z_2)+
 2t_1 t_2g_B(z_1,
z_2, \gamma)\right] \left(\frac{1}{\log
r}\right)^2+o\left(\left(\frac{1}{\log r}\right)^2\right),
\]
\[
 r \to 0.
\]
We applied here the symmetry of the Robin function:\,\,\, $g_B(z_1, z_2, \gamma)=$ $g_B$
$(z_2, z_1,\gamma)$. Again by formula (1.1) we have
\[
\mathrm{cap}\, C(r; G,\Gamma, W, \Delta, \Psi)=2\pi(t_1^2+t_2^2)
\left(-\frac{1}{\log r}\right)-2\pi\left[t_1^2\log \frac{r(G,
\Gamma, w_1)}{|f'(z_1)|}+\right.
\]
\[
\left.+t^2_2\log \frac{r(G, \Gamma, w_2)}{|f'(z_2)|}+2t_1t_2
g_G(w_1, w_2, \Gamma)\right]\left(\frac{1}{\log
r}\right)^2+o\left(\left(\frac{1}{\log r}\right)^2\right),
\]
$
r \to 0.
$
Summing up the above relations for the capacities of condensers 
we arrive at the inequality (2.1).

We now assume, that the inclusion 
$f((\partial B)\backslash\gamma)\subset(\partial G)\backslash$ $\Gamma$ 
holds, and let the function $u$ be continuous at the points of 
$\overline{f(B)}$, and satisfy Lipschitz condition in a neighborhood 
of every finite point of the domain $f(B),$ equal to zero in 
a neighborhood of $f(\gamma)$ and equal to $t_k$ in a neighborhood of 
the sets $E(w_k, |f'(z_k)|r),\, k=1, 2$. The function $\widetilde{u}$, 
defined as an extension of $u$:
\[
\widetilde{u}(w) =\left\{
\begin{array}{lll}
u(w)\,\,\,\,,\quad w\in\overline{f(B)} \\
\quad 0\,\,\,\,\,\,\,, \quad w\in\overline{G}\backslash \overline{f(B)},   \\
 \end{array}
 \right.
\]
satisfies the aforementioned properties with $f(B)$ in place of 
$G$ and $f(\gamma)$ in place of $\Gamma.$ Therefore
\[
I(u, f(B))= I(\widetilde{u}, G)\geq \mathrm{cap}\, C(r; G,\Gamma,
W, \Delta, \Psi).
\]

In view of the arbitrary choice of the function $u$, we have
\[
\mathrm{cap}\, C(r; f(B), f(\gamma),  W, \Delta, \Psi)\geq
\mathrm{cap}\, C(r; G,\Gamma, W, \Delta, \Psi).
\]
Therefore, unlike in the first case, now also an inequality in 
the other direction is valid. Repeating again a corresponding 
part of the proof of inequality (2.1), we arrive at the 
inequality (2.2). The theorem is proved.

 \vspace{0.5cm}{\bf 2.2. Remark.} From the result of the paper [DubPr] 
it follows that the sign of the equality in (2.1) is attained only in 
the case when $\overline{f(B)}=\overline{G},$ $f(\gamma)=\Gamma$ and 
the set $G\cap\partial f(B)$ consists of a finite number of piecewise 
smooth curves, and at the interior points of the curves the normal 
derivative of the function $t_1g_G(w, w_1, \Gamma)+ t_2g_G(w, w_2, \Gamma)$ 
is equal to zero. The equality in (2.2) hold if and only if 
$\overline{f(B)}=\overline{G}$ and $t_1g_G(w, w_1, \Gamma)+ t_2g_G(w, w_2, \Gamma)=0$ 
on $G\cap f(\gamma)$.


 \vspace{2.0cm}
  \noindent {\Large\bf 3 Particular cases}

\vspace{0.5cm} First of all we observe that for $t_1=1,\, t_2=0$ 
Theorem 2.1 has the character of a majorization principle: 
if $f(\gamma)\subset\Gamma$, then for every point $z$ in 
the domain $B$ the following holds
\[
\qquad\qquad\qquad \qquad\qquad\qquad |Df(z)|\geq 1,
\qquad\qquad\qquad\qquad\qquad\qquad (3.1)
\]
whereas in the case 
$f((\partial B)\backslash\gamma)\subset(\partial G)\backslash\Gamma$ we have
\[
\qquad\qquad\qquad \qquad\qquad\qquad |Df(z)|\leq 1,
\qquad\qquad\qquad\qquad\qquad\qquad (3.2)
\]
for every point $z \in B$. In particular, if $\gamma=\partial B,$ 
$\Gamma=\partial G$ the inequality (3.2) is well-known as the 
monotonicity property of the inner radius of the domain. 
If $\gamma\not=\partial B,$ $\Gamma\not=\partial G$, 
then the inequalities (3.1) and (3.2) can be extended to 
a boundary point $z\in (\partial B)\backslash\gamma$ under 
the condition of the existence of the corresponding limit 
of the derivative $f'$. Let us now consider the particular 
cases of the inequalities (2.1), (2.2), (3.1), (3.2), 
when $B=\{z:|z|<1\}$ and $G=\{w:|w|<1\}.$

\vspace{0.5cm} {\bf  3.1. Corollary.}  {\it Let the function} $f$ 
{\it be regular and univalent in the disk} $|z|<1,$ $f(0)=0$ 
{\it and} $|f(z)|<1,$ {\it for} $|z|<1.$ {\it We suppose that the function} 
$f$ {\it maps in the sense of boundary correspondence a set} $\alpha$, 
{\it consisting of a finite number of open arcs of the circle} $|z|=1,$ 
{\it into the circle} $|w|=1.$ {\it Then}
\[
\sqrt{|f'(0)|}\,\mathrm{cap}\, \{w:|w|=1, \,\, w\notin
f(\alpha)\}\leq \mathrm{cap}\,\{z:|z|=1, \,\, z\notin \alpha\},
\]
 {\it where} $\mathrm{cap}\,(\cdot)$ {\it denotes the logarithmic capacity}.

\vspace{0.5cm} {\bf  Proof.} We may assume that the sets 
$\gamma=\{z:|z|=1,\, z\notin \alpha\}$ and $\Gamma=\{w:|w|=1, \,\,
w\notin f(\alpha)\}$ are admissible. The inequality (3.2) gives
\[
r(\{z:|z|<1\},\gamma, 0)|f'(0)|\leq r(\{w:|w|<1\}, \Gamma, 0).
\]
It remains to verify that
 \[
\qquad\qquad\qquad \qquad r(\{z:|z|<1\},\gamma, 0)=(
\mathrm{cap}\,\gamma)^{-2}. \qquad\qquad\qquad\qquad
 (3.3)
 \]
In view of the symmetry with respect to the circle $|z|=1$ we have
\[
g_U(z, 0,\gamma)=g_{\overline{\mathbf{C}}_z\backslash\gamma}(z,
0)+g_{\overline{\mathbf{C}}_z\backslash\gamma}(z, \infty),
\]
where $U=\{z:|z|<1\}$. Therefore
\[
\log r(U, \gamma, 0)= \log r
(\overline{\mathbf{C}}_z\backslash\gamma, 0)+
g_{\overline{\mathbf{C}}_z\backslash\gamma}(0, \infty).
\]
On the other hand
\[
g_{\overline{\mathbf{C}}_z\backslash\gamma}(z, 0)+ \log |z|=
g_{\overline{\mathbf{C}}_z\backslash\gamma}(z, \infty)
\]
and therefore
\[
\qquad\qquad\qquad\qquad\log
r(\overline{\mathbf{C}}_z\backslash\gamma, 0)=
g_{\overline{\mathbf{C}}_z\backslash\gamma}(0, \infty).
\qquad\qquad\qquad\qquad\quad(3.4)
\]
Summing up the above equalities we get
\[
r(U,\gamma, 0)=(r(\overline{\mathbf{C}}_z\backslash\gamma, 0))^2,
\]
which is equivalent to (3.3). The corollary is proved.

Applying the inequality (3.1) instead of (3.2) we arrive as above to 
the following result of Pommerenke.

\vspace{0.5cm} {\bf  3.2. Corollary} [P, p. 217].  {\it If the function} 
$f$ {\it is regular and univalent in the disk} $|z|<1,$ $f(0)=0$ 
{\it and} $|f(z)|<1$ {\it for} $|z|<1,$ {\it and if a closed set}
 $\gamma$ {\it on the circle} $|z|=1$ 
{\it is mapped in the sense of boundary correspondence under the mapping} 
$f$ {\it onto a closed set} $\Gamma$ {\it on the circle} $|w|=1,$
{\it then}
\[
\sqrt{|f'(0)|}\,\mathrm{cap}\,\Gamma \geq \mathrm{cap}\,\gamma.
\]

\vspace{0.5cm} {\bf  3.3. Corollary} [Neh].  {\it Let the function}
$f$ {\it map the disk} $|z|<1$ {\it conformally and univalently into the disk}
 $|w|<1.$ {\it Then for all points} $z_1, z_2$ {\it of the disk} $|z|<1$ 
{\it and for all real numbers} $t_1, t_2$ 
{\it the following inequality is valid}
\[
\left[\frac{(1-|z_1|^2)|f'(z_1)|}{1-|f(z_1)|^2}\right]^{t_1^2}\,\,
\left[\frac{(1-|z_2|^2)|f'(z_2)|}{1-|f(z_2)|^2}\right]^{t_2^2}\leq
\]

\[
\leq\left|\frac{(z_1-z_2)(1-{\overline{f(z_1)}}
f(z_2))}{(1-{\overline z_1}z_2)(f(z_1)-f(z_2))}\right|^{2t_1t_2}.
\]
This inequality follows immediately from (2.2). 
For some applications of this inequality see [KrR, p. 125]. 
For $t_1=1, t_2=0$ we obtain the inequality of Pick
\[
|f'(z)|(1-|z|^2)\leq  1-|f(z)|^2, \,\,\,\, |z|<1,
\]
which, as well-known, is valid for all funtions $w=f(z),$ $|f(z)|<1$, 
for $|z|<1$ regular in the disk $|z|<1$ (not necessarily univalent). 
Setting in Corollary 3.3 $z_1=0,$ $z_2= re^{i\varphi},$ $t_1=-t_2=1,$ 
and letting $r \to 0,$ we arrive after simple calculations at the 
following statement for the Schwarzian derivatives $S_f(z).$

\vspace{0.5cm} {\bf  3.4. Corollary}.  {\it If the function} 
$f(z)=c_1z+\ldots$ {\it is regular and univalent in the disk} 
$|z|<1$ {\it and} $|f(z)|<1$ {\it for} $|z|<1,$ {\it then}
\[
|S_f(0)|\leq 6(1-|c_1|^2).
\]

 \vspace{0.5cm} {\bf  3.5. Corollary}. {\it Let the function} $f$ 
{\it be regular and univalent in the disk} $|z|<1,$ {\it  let}$|f(z)|<1$ 
{\it for} $|z|<1$, {\it and let} $f(\gamma)\subset\Gamma$, 
{\it where} $\gamma$ {\it is a closed subset of the circle} 
$|z|=1,$ {\it different from the whole circle} $|z|=1,$ {\it and let} 
$\Gamma$  {\it be a similar subset of
} $|w|=1.$ {\it Then for all } $z_1, z_2$ {\it in}
 $|z|<1$ {\it and all real numbers} $t_1, t_2$ {\it the following
inequality holds }

\[
\prod\limits_{k=1}^2\left[\frac{r^2(\overline{\mathbf{C}}_z\backslash\gamma,
z_k)(1-|f(z_k)|^2)|f'(z_k)|}
{r^2(\overline{\mathbf{C}}_w\backslash\Gamma,
f(z_k))(1-|z_k|^2)}\right]^{t^2_k} \geq
\]

\[
\geq
\left\{\frac{\exp[g_{\overline{\mathbf{C}}_{w\backslash\Gamma}}(f(z_1),
f(z_2))+ g_{\overline{\mathbf{C}}_{w\backslash\Gamma}}(f(z_1),1/
\overline{f(z_2)}\,)]}
{\exp[g_{\overline{\mathbf{C}}_{z\backslash\gamma}}(z_1,
z_2)+g_{\overline{\mathbf{C}}_{z\backslash\gamma}}(z_1,
1/{\overline{z}_2)]}} \right\}^{2t_1t_2}.
\]

 \vspace{0.5cm} {\bf  Proof.} According to Theorem 2.1
\[
\prod\limits_{k=1}^2\left[\frac{r(U_z, \gamma, z_k)|f'(z_k)|}
{r(U_w,\Gamma, f(z_k))}\right]^{t_k^2} \geq
\left\{\frac{\exp[g_{U_w}(f(z_1), f(z_2),
\Gamma)]}{\exp[g_{U_z}(z_1, z_2, \gamma)]}\right\}^{2t_1t_2},
\]
where $U_z=\{z:|z|<1\}$ and $U_w=\{w:|w|<1\}$. Applying the inequality
\[
g_{U_z}(\zeta, z,
\gamma)=g_{\overline{{\mathbf{C}}_z}\backslash\gamma}(\zeta, z)+
g_{\overline{{\mathbf{C}}_z}\backslash\gamma}(\zeta,1/\overline{z}),
\,\, \zeta, z \in U_z,
\]
we conclude that the right hand side of the above inequality agrees 
with the right hand side of the inequality of Corollary 3.5. 
Furthermore, from this we conclude that
\[
\log r(U_z, \gamma, z_k)=\log r
(\overline{\mathbf{C}}_{z}\backslash\gamma, z_k)+
g_{\overline{\mathbf{C}}_{z}\backslash\gamma}(z_k,
1/{\overline{z}}_k), \,\,\, k=1, 2.
\]
Let the function $\varphi(z):=(z-z_k)/(1-{\overline{z}}_k z)$ and 
$B:=\varphi(\overline{\mathbf{C}}_{z}\backslash\gamma)$. 
From the conformal invariance of the Green function and the 
inequality (3.4) it follows that
\[
g_{\overline{\mathbf{C}_z}\backslash\gamma} (z_k,
1/{\overline{z}}_k)= g_B(0, \infty)= \log r (B, 0)=
\]

\[
 \log [r(\overline{\mathbf{C}}_{z}\backslash\gamma,
 z_k)|\varphi'(z_k)|\,]= \log[r(\overline{\mathbf{C}}_{z}\backslash\gamma, z_k)/
 (1-|z_k|^2)\,],\,\,\, k=1,2.
 \]
Therefore
\[
r(U_z, \gamma, z_k)=
\frac{r^2(\overline{\mathbf{C}}_{z}\backslash\gamma,
z_k)}{1-|z_k|^2}, \,\,\, k=1, 2.
\]
Similar representation holds for $r(U_w, \Gamma, f(z_k))$,  $k=1, 2,$ 
which completes the proof of Corollary 3.5.

 \vspace{0.5cm} {\bf  3.6. Corollary}. {\it Under the assumptions of
 Corollary} 3.5 {\it we assume furthermore that at two points} $z_k, |z_k|=1,$
$z_k\notin\gamma$ {\it the angular limits} $f(z_k)$ {\it exist with} $|f(z_k)|=1$.
{\it Then for the  angular derivatives} $f'(z_k)$ {\it the following inequality holds}
\[
\prod\limits_{k=1}^2\left[\frac{r(\overline{\mathbf{C}}_z\backslash\gamma,
z_k)} {r(\overline{\mathbf{C}}_w\backslash\Gamma,
f(z_k))}|f'(z_k)|\right]^{t^2_k} \geq
\exp\left\{2t_1t_2[g_{\overline{\mathbf{C}}_{w\backslash\Gamma}}(f(z_1),
f(z_2))-\right.
\]

\[{\hskip-10cm\left.
-g_{\overline{\mathbf{C}}_{z\backslash\gamma}}(z_1, z_2)\,\,
]\right\}}
\]
{\it for all real} $t_1$ {\it and}  $t_2$.

The proof follows from Corollary 3.5 with a limiting passage of points in 
$|z|<1$ to boundary points (cf. [P, p. 79--83]). We next represent a corollary 
of Theorem 2.1 for functions defined in an annulus.

 \vspace{0.5cm} {\bf  3.7. Corollary}. {\it Let the function} $f$ 
{\it be regular and univalent in the annulus} $R=\{z:\rho<|z|<1\},$ 
{\it whose image} $f(R)$ {\it lies in the disk} $|w|<1$ {\it and} 
$|f(z)|=1$ {\it  for} $|z|=1$ {\it and let} $\beta$ 
{\it be a closed subset of the circle} $|z|=1,$ 
{\it consisting of a finite number of nondegenerate arcs.} 
{\it Then for all points} $z_1$ {\it and} $z_2,$ {\it in the set} 
$R\cup\{z:|z|=1, z \notin\beta\}, $ {\it the following inequality holds}
\[
\prod\limits_{k=1}^2\frac{r(R, \gamma, z_k)|f'(z_k)|} {r(U_w,
f(\beta), f(z_k))}
 \leq
\exp\left\{ 2[g_R((z_1, z_2, \gamma)-g_{U_w}(f(z_1), f(z_2),
f(\beta))\, ]\right\},
\]
{\it where} $\gamma=\beta\cup\{z:|z|=\rho\},$ $U_w=\{w:|w|<1\}.$

\vspace{0.5cm} {\bf  Proof.} Setting in Theorem 2.1 $B=R$, $G=U_w,$ 
$\Gamma=f(\beta)$ and $t_1=-t_2=1,$ we obtain from (2.2) the desired relation.

The particular case of Corollary 3.7, when $\beta=\emptyset,$ $|z_1|=|z_2|=1$ 
was considered by A.Yu. Solynin [Sol2, p. 135]. For $\beta=\{z:|z|=1\}$ 
this corollary coincides with Theorem 1.4 of [DubKos], which in the methodic 
sense goes back to Nehari [Neh]. Distortion theorems for functions 
described in Corollary 3.7 were also earlier considered in [D2, Huck].

 \vspace{0.5cm} {\bf  3.8. Corollary} [Sin, A].
 {\it  Under the hypotheses of Corollary} 3.7 {\it let the set} 
$\beta$ {\it coincide with the circle} $|z|=1.$
 {\it Then for every point } $z$ {\it of the ring} $R$ 
{\it we have the inequality}
 \[
 \left|\frac{1}{4}S_f(z)+\pi l(z, z) \right|\leq \pi K(z,
 \overline{z})-\frac{|f'(z)|^2}{(1-|f(z)|^2)^2},
 \]
 {\it where} $K (\cdot, \cdot)$ {\it and} $l (\cdot, \cdot)$ 
{\it are Bergman kernels of} 
 {\it the first and second kind of the domain} $R$ 
{\it with respect to the class of
all functions, regular  in }$R$ {\it and with square
integrable modulus in } $R.$ 

 \vspace{0.5cm} {\bf  Proof.} Let $z_0$ be an arbitrary fixed point 
of the annulus $R$ and let $\varphi$ be a real number. 
Applying Corollary 3.7 in the case of the points
 $z_1=z_0+\rho e^{i\varphi}$ and $z_2=z_0-\rho e^{i\varphi}$, we get

 \[
 \frac{r(R, z_0+\rho e^{i\varphi})r(R, z_0-\rho e^{i\varphi})}
{4\rho^2 \exp[2g_R(z_0+\rho e^{i\varphi},z_0-\rho e^{i\varphi})]}\cdot
 \frac{|f'(z_0+\rho e^{i\varphi})f'(z_0-\rho e^{i\varphi})|4 \rho^2}
 {|f(z_0+\rho e^{i\varphi})-f(z_0-\rho e^{i\varphi})|^2}\leq
 \]

 \[
\qquad\qquad\qquad \leq\frac{(1-|f(z_0+\rho
e^{i\varphi})|^2)(1-|f(z_0-\rho e^{i\varphi})|^2)}{|1-f(z_0+\rho
e^{i\varphi}) \overline{f(z_0-\rho e^{i\varphi}})|^2}.
 \qquad\qquad(3.5)
\]
Making use of well-known relations between Bergman kernels 
and Green functions (see, e.g.
[S]), as well as the Taylor expansion of the function $f$ 

\[
f(z)= c_0+c_1(z-z_0)+c_2(z-z_0)^2+c_3(z-z_0)^3+ \ldots,
\]
we conclude that the first quotient of the left side of the
inequality (3.5) is equal to
 \[
 1-4\pi\{K(z_0, \overline{z}_0)-\mathrm{Re}[e^{2i\varphi} l(z_0,
 z_0)]\}\rho^2+o(\rho^2), \,\,\, \rho \to 0.
 \]
 The second quotient of the left side of (3.5) has the form

 \[
 \left|1+6\left(\frac{c_3}{c_1}-\frac{c_2^2}{c_1^2}
 \right)e^{2i\varphi}\rho^2+o(\rho^2)
 \right|=
 \left|1+S_f(z_0)e^{2i\varphi}\rho^2+o(\rho^2)
 \right|,
 \,\,\, \rho \to 0.
 \]
Finally the right side of (3.5) after simple transformations 
can be expressed in the form 
 \[
 1-\frac{4|c_1|^2}{(1-|c_0|^2)^2}\,\rho^2+o(\rho^2), \,\,\, \rho \to
 0.
 \]
Adding together the above expressions we arrive at the 
required inequality  for $z=z_0$ because $\varphi$ was arbitrary.





\newpage
\noindent {\Large\bf 4 Majorization principles for
regular functions}

\vspace{0.5cm}\noindent We begin with the analogue of Lindel\"of's principle
which expresses the behavior of the Robin function under 
a regular mapping. The proof of this statement resembles
in many respects the proof of the Lindel\"of principle itself for the Green
function
[Sto,Ch. Y, \S2].

\vspace{0.5cm} {\bf  4.1. Theorem.}  {\it Let the domains} $B$ {\it
and} $G$, {\it lying in the planes} $\mathbf{C}_z$
{\it and} $\mathbf{C}_w$, {\it respectively, and let also the subsets}
$\gamma\subset\partial B$ {\it and} $\Gamma\subset\partial G$ 
{\it be admissible. We assume that the function} $f$ {\it is regular in the
domain}
$B,$ $f(B) \subset G$ {\it and} $f((\partial
B)\backslash\gamma)\subset(\partial G)\backslash\Gamma$ ({\it 
i.e for each sequence of points} $\zeta_n \in B,$ {\it
approaching the set
} $(\partial B)\backslash\gamma,$ {\it
the corresponding sequence} $f(\zeta_n) \to
(\partial G)\backslash\Gamma$). {\it Let} $w_0$ {\it be a point of} $f(B),$ 
let $\{z_\nu\}$ $(\nu=0, 1, \ldots)$ {\it be in } $B,$ {\it with} $f(z_\nu)=w_0$ 
{\it and} $n_\nu$ let
{\it be the orders at the points} $z_\nu$ {\it of the zeros of} $f(z)-w_0.$ {\it Then}

\[
\qquad\qquad\qquad\qquad g_G(f(z), w_0, \Gamma) \geq
\sum\limits_{\nu \geq 0}\, n_\nu g_B(z, z_\nu, \gamma)
 \qquad\qquad\qquad(4.1)
\]
{\it for all} $z \in B.$ {\it Equality in} (4.1) {\it for a single point } 
$z \in B$ {\it implies that}
(4.1) {\it holds at every point} $B$.

 \vspace{0.5cm} {\bf  Proof.} Applying, if necessary, a conformal and univalent mapping 
we may suppose without restriction of generality
that the domains $B$  and $G$  are bounded by a finite number of
circles and that the sets   
$\gamma$ and $\Gamma$ consist of a finite number of arcs on these
circles. In this case the function 
$f$ is defined  on $B\cup(\partial B)\backslash \gamma$.
We fix a natural number $N,$ and consider the function
  \[
  I_{N}(z)= g_G(f(z), w_0, \Gamma)-
  \sum\limits_{\nu=0}^{N}\, n_\nu g_B(z, z_\nu,\gamma),
  \]
  defined on the set $B\backslash\bigcup\limits_{\nu\geq 0}z_\nu$.
  Applying the expansion of the function $f$ in a neighborhood of the
points $z_\nu$:

  \[
  f(z)= w_0+c_{n_\nu}(z-z_\nu)^{n_\nu}+\ldots\qquad,
  c_{n_{\nu}}\not= 0,
  \]
  and also the representation of the Robin function in a neighborhood of
a pole, we easily conclude that the points 
$z_\nu$, $\nu \leq N,$ 
are removable singularities of the function $I_{N}(z)$. This function
$I_{N}(z)$ is harmonic in the domain
$B\backslash\bigcup\limits_{\nu>N}z_\nu,$ approaches $+\infty$ when
$z \to z_\nu$, $\nu >N,$ is nonnegative on the boundary points at $\gamma$, and at the
points  $z \in (\partial B)\backslash\gamma$ satisfies the condition:
$\partial I_{N}(z)/\partial n=0.$
   By the maximum
  principle of E. Hopf we conclude that $I_{N}(z) \geq 0$ in the domain 
$B\backslash\bigcup\limits_{\nu>N}z_\nu,$
  (cf., for instance, [ProW]).
 In view of the arbitrariness of $N$, the inequality (4.1) holds for all
 $z \in B.$ At the same time we proved that the series
 $\sum\limits_{\nu\geq 0}n_{\nu}  g_B(z, z_\nu, \gamma)$ converges and that the
function
\[
\qquad\qquad\qquad\qquad g_G(f(z), w_0, \Gamma) - \sum\limits_{\nu
\geq 0}\, n_\nu g_B(z, z_\nu, \gamma)
 \qquad\qquad\qquad(4.2)
\]
is nonnegative and harmonic
 in $B$. By virtue of the maximum principle if the function is zero at a point
of the domain $B$ then it vanishes identically.
The theorem is proved.

In the papers [M1], [M2]  Mityuk introduced into consideration 
a theorem about the change
of the interior radius of the domain under a regular mapping and proved the
effectiveness of this result with symmetrization methods.
Following Mityuk, we consider the corresponding majorazation principle for
the quantity
 $r(B, \gamma, z_0)$, which also may be considered as a distortion theorem.

  \vspace{0.5cm} {\bf  4.2. Theorem.}  {\it Under the hypotheses of Theorem}
4.1 {\it suppose that in a neighborhood of the point
} $z_0$ {\it we have the expansion}
\[
f(z)= w_0+c_{n}(z-z_0)^{n}+\ldots\qquad,
  c_{n}\not= 0,
  \]
$(n=n_0).$ {\it Then}

\[
\quad\quad r(G, \Gamma, w_0) \geq |c_n| r^n(B, \gamma, z_0)
\exp\left\{\sum\limits_{\nu\geq 1} n_\nu g_B(z_0, z_\nu, \gamma)
\right\}.
 \quad\quad\quad(4.3)
\]
{\it If the mapping } $f$ {\it satisfies }
$f(\gamma)=\Gamma$ {\it and} $f((\partial
B)\backslash\gamma)=(\partial G)\backslash\Gamma,$  {\it then equality holds}
{\it in the formula} (4.3). 

 \vspace{0.5cm} {\bf  Proof.} From inequality (4.1)
it follows that in a small enough neighborhood of
$z_0$

 \[
 -\log|c_n(z-z_0)^n +\ldots|+ \log r(G, \Gamma, w_0)+ o(1) \geq -
 n \log|z-z_0|+
 \]
 \[
 +n \log r(B, \gamma, z_0) +\sum\limits_{\nu\geq 1}n_{\nu}  g_B(z, z_\nu,
 \gamma), \qquad z \to z_0.
 \]
 Therefore
 \[
 \log r(G, \Gamma, w_0)\geq \log|c_n| +\log|1+o(1)|+\log r^n(B,
 \gamma, z_0)+
 \]
 \[
  +\sum\limits_{\nu\geq 1}n_{\nu}  g_B(z_0, z_\nu,
 \gamma)+o(1), \qquad z \to z_0.
 \]
Passing to the limit when
$z \to z_0,$ we obtain the inequality (4.3). If $f(\gamma)=\Gamma$ and
 $f((\partial B)\backslash\gamma)=(\partial G)\backslash\Gamma$,
 then by the maximum principle of E. Hopf we conclude that 
the function (4.2) vanishes identically.
Consequently, in (4.1) and also in (4.3)
the equality sign holds. The theorem is proved.

 From Theorems 4.1 and 4.2 
there follows a statement concerning the behavior
of the quadratic form
 \[
\sum\limits_{k=1}^n\, t^2_k \log r(B, \gamma, z_k)+
\sum\limits_{k=1}^n\sum\limits_{l=1\atop l\not=k}^n\, t_k t_l
g_B(z_k, z_l, \gamma)
\]
under a regular mapping $f$ in the case when the points $z_k$ of the set 
$Z=\{z_k\}_{k=1}^n$
 are located in the domain $B$, and the real numbers $t_k$ of the set 
$\Delta=\{t_k\}_{k=1}^n$, have the same sign. For arbitrary real numbers
 $t_k$ and $p$-valent
 functions $f$
the following counterpart of the majorization principle of
 [Dub2] holds.

  \vspace{0.5cm} {\bf  4.3. Theorem.}  {\it We suppose that in the
hypotheses of Theorem} 4.1
{\it the function} $w=f(z)$ is $p$-{\it valent in} $B$. {\it Let} 
$w_l,l=1, \ldots, m,$  {\it be distinct points of the domain } $f(B)$, {\it
each of which has exactly} $p$ {\it preimages in} $B$ {\it taking into account multiplicity,
and let } 
$z_{jl},$ $j=1, \ldots, p,$ {\it
be preimages of the points} $w_l, l=1, \ldots, m$ {\it (each zero of the function
} $f(z)-w_l$ {\it occurs as many times as its multiplicity indicates).
 }.{\it Then for all real numbers
 } $t_l, l=1,
\ldots, m,$ {\it the following inequality holds}

\[
\qquad\qquad p\left\{\sum\, t^2_l \log r(G, \Gamma, w_l)+\sum
t_kt_l g_G(w_k, w_l, \Gamma) \right\}\geq
\]

\[
\geq \sum t_l^2[\log |c_{jl}|+ p_{jl}\log r(B, \gamma,
z_{jl})]+\sum t_kt_l g_B(z_{ik}, z_{jl}, \gamma), \qquad(4.4)
\]
{\it where} $c_{jl}$ {\it and} $p_{jl}$ {\it are defined from the expansion}

\[
f(z)-w_l= c_{jl}(z-z_{jl})^{p_{jl}}+\ldots,\,\,\, c_{jl}\not=0,
\]
$j=1, \ldots, p, \, l=1, \ldots, m$ ({\it here and later the symbol}
$\sum$ {\it denotes summation,
over all indices appearing from the context, excluding those for which the summand
is either } $\infty$ {\it or not defined}). {\it If, furthermore,
} $f(\gamma)=\Gamma,$ $f((\partial
B)\backslash\gamma)=(\partial G)\backslash\Gamma$ {\it 
and the mapping} $f$ {\it defines a} $p$-{\it valent
covering of the domain} $G$ {\it onto} $B$, {\it then in the inequality}
(4.4) {\it the equality sign holds.}

\vspace{0.5cm} {\bf Proof.} 
 We may assume that the domains $B$ and  $G$ are bounded by
a finite number of circles and all numbers $t_l,$ $l=1, \ldots, m,$ are nonzero.
Set $W=\{w_l\}$,
$\Delta=\{t_l\},$ $\Psi=\{\psi_l(r)\},$ $\psi_l(r)\equiv r, l=1,
\ldots, m,$ and consider the condenser
$C(r; G, \Gamma, W, \Delta,
\Psi)$. Our task is to compare the capacity of this condenser with 
the capacity of the condenser $\widetilde{C}(r)$,
 defined as follows.
Let $\widetilde{Z}=\{\widetilde{z}_{kl}\},$ $\widetilde{z}_{kl},
\, k=1, \ldots, n_l,$ $l=1, \ldots, m,$ { be the zeros of the function}$f(z)-w_l$
 { without counting the multiplicity
} $\widetilde{\Delta}=\{t_{kl}\},$ $t_{kl}=t_l / p, \,
k=1, \ldots, n_l,$; $\widetilde{\Psi}=\{\psi_{kl}(r)\},$
$\psi_{kl}(r)=|\widetilde{c}_{kl}|^{-1/p_{kl}} r^{1/p_{kl}},\,
$ and let $p_{kl}, k=1, \ldots, n_l,$ 
be the multiplicities at the points $\widetilde{z}_{kl},$
$\sum^{n_l}_{k=1}\,$ $p_{kl}=p, \, l=1, \ldots, m,$
$\widetilde{c}_{kl}$ be the coefficient, corresponding to the point
$\widetilde{z}_{kl}$. The condenser $\widetilde{C}(r)$ is obtained
from the condenser $C(r; B, \gamma, \widetilde{Z}, \widetilde{\Delta},
\widetilde{\Psi})$ by changing the disks $E(\widetilde{z}_{kl},
\psi_{kl}(r))$ with almost disks $\widetilde{E}(\widetilde{z}_{kl},
\psi_{kl}(r))$. Here $\widetilde{E}(\widetilde{z}_{kl},
\psi_{kl}(r))$ is the connected part of the preimage of 
$E(w_l, r)$ under the mapping 
$f$, lying in the neighborhood of the point $\widetilde{z}_{kl}$.
Let now  $u$ be the potential function of the condenser
$\widetilde{C}(r)$,
i.e. the real-valued function $u$, continuous in $\overline{B}$, 
harmonic in $B\backslash  \bigcup_{k,l}$
$\widetilde{E}(\widetilde{z}_{kl},\psi_{kl}(r))$, equal to $0$ on
$\gamma$ and equal to $t_{kl}$ 
on $\widetilde{E}(\widetilde{z}_{kl},\psi_{kl}(r))$ 
and satisfying the condition $\partial u/\partial n=0$  at the
points of $(\partial
B)\backslash\gamma$. On the set $f(B)$ we define the function 

\[
U(w)=\sum\limits_{f(z)=w}\, u(z).
\]
For each
$w \in f(B)$ the value of the function $U(w)$ is the sum
of  
at most $p$ summands. From the definition of the capacity of a condenser, the
convexity of the function  $y=x^2$ 
and conformal invariance of the Dirichlet integral
we have
\[
\mathrm{cap}\,C(r; G, \Gamma, W, \Delta, \Psi) \leq \iint_{f(B)}\,
|\nabla U|^2 dudv\leq p \iint_{B} |\nabla u|^2 dxdy.
\]
Applying the Dirichlet principle we conclude that

\[
\mathrm{cap}\,C(r; G, \Gamma, W, \Delta, \Psi) \leq p\,\,
\mathrm{cap}\,\widetilde{C}(r).
\]

Applying the asymptotic formula (1.1)  to both condensers we have

\[
p\left\{\sum t^2_l\log r(G, \Gamma, w_l)+\sum t_k t_l g_G(w_k,
w_l, \Gamma)\right\}\geq
\]

\[
\geq \sum t^2_l [p_{kl}\log|\widetilde{c}_{kl}|+p_{kl}^2 \log r(B,
\gamma, \widetilde{z}_{kl})]+ \sum t_k  p_{ik}t_l
p_{jl}g_B(\widetilde{z}_{ik}, \widetilde{z}_{jl}, \gamma) ,
\]
which coincides with the inequality (4.4).

Let now  the function $f$ define  a complete $p$-to-one covering mapping of the
domain
$G$ onto $B$ and $f(\gamma)=\Gamma,$ $f((\partial
B)\backslash\gamma)=(\partial G)\backslash\Gamma$. If
$\omega(w)$ is the potential function of the condenser $C(r; G, \Gamma,
W, \Delta, \Psi)$, then the composite function $\omega(f(z))$ is the
potential for the condenser $\widetilde{C}(r)$, 
and by the conformal invariance of the
Dirichlet integral and the   $p$-valence of the covering mapping we have

\[
\mathrm{cap}\,C(r; G, \Gamma, W, \Delta, \Psi)= p\,\,
\mathrm{cap}\,\widetilde{C}(r).
\]

Taking into account the formula
(1.1) , this gives the equality sign in (4.4). The theorem is
proved. It can be proved that in the case when the
$t_l$'s have the different signs, the
$p$-valence in Theorem 4.3 is essential ([Dub2, p.
538]).


 \vspace{1.5cm}\noindent {\Large\bf 5  Examples}

\vspace{0.5cm}\noindent As an application of the majorization principles of the
preceding section we now consider some corollaries to
 the
distortion in regular mappings. Immediately from Theorem 4.1 it follows that

\vspace{0.5cm} {\bf  5.1. Corollary.}  {\it Let the domain}
$B\subset\mathbf{C}_z$ {\it and} $G\subset\mathbf{C}_w$,  {\it and also the sets } 
$\gamma\subset\partial B,$ {\it }
$\Gamma\subset\partial G$ {\it be admissible. We require that the function} 
$f$ {\it is regular in the domain} $B,$ $f(B) \subset G$ {\it
and} $f((\partial B)\backslash\gamma)\subset(\partial
G)\backslash\Gamma.$  {\it Then for every pair of distinct points}
$z$ {\it and} $\zeta$ {\it of the domain} $B$ {\it the following inequality is valid}
\[
\qquad\qquad\qquad\qquad \delta(f(z), f(\zeta), G,  \Gamma) \leq
\delta(z, \zeta, B,  \gamma).
 \qquad\qquad\qquad(5.1)
\]

{\it The equality in }(5.1) {\it is attained in the case of conformal
and univalent functions} $f$.

 If the domains $B$  and $G$ are disks, and the sets  $\gamma$  and $\Gamma$ 
are the boundaries of these disks, then  (5.1) coincides with the invariant form of the  
Schwarz lemma 
due to G. Pick. We give an example of the inequality in (5.1) for the case when
$B=\{z:0<\mathrm{Im}z< \pi/2\},$ $\gamma=\{z:\mathrm{Im}z=
\pi/2\},$ $G=\{w:0<\mathrm{Im}w< \pi/2\}$ and
$\Gamma=\{w:\mathrm{Im}w= \pi/2\}$:

\[
\left|\frac{\big(e^{f(z)}-e^{f(\zeta)}\big)\big(e^{f(z)}-e^{\overline{f(\zeta)}}\big)}
{\big(e^{f(z)}+e^{f(\zeta)}\big)\big(e^{f(z)}+e^{\overline{f(\zeta)}}\big)}
\right|\leq \left|\frac{(e^z-
e^\zeta)(e^z-e^{\overline{\zeta}})}{(e^z+
e^\zeta)(e^z+e^{\overline{\zeta}})} \right|,
\]
for all $z, \zeta,$ $0<\mathrm{Im}z< \pi/2$,
$0<\mathrm{Im}\zeta< \pi/2$.

A particular case of  Theorem  4.2 is the following statement.

\vspace{0.5cm} {\bf  5.2. Corollary.}  {\it In the hypotheses of Corollary}
5.1 {\it let the point} $z_0 \in B,$ $w_0=f(z_0);$ {\it and let} $z_\nu,\, \nu=1,
2, \ldots,${\it be the zeros of the function} $f(z)-w_0,$ {\it different from}
$z_0$ {\it and let} $n_\nu$ {\it be the multiplicities of}
 $z_\nu,\, \nu=1, 2,
\ldots.$ {\it Then}
\[
\qquad\qquad\quad   |D f(z_0)|\leq \exp
\left\{-\sum\limits_{\nu\geq 1}\, n_\nu\, g_B(z_0, z_{\nu},
\gamma) \right\} \leq 1.
 \qquad\qquad (5.2)
\]

Therefore, the inequality (3.2) holds for arbitrary regular functions
(not necessarily univalent) and even allows a refined formulation 
taking into account the multiplicity.
In particular, if the function $f$ is regular in the unit disk $|z|<1$, $f(0)=0,$ 
and $|f(z)|<1$ for
$|z|<1$, and if $f$ maps the set $\alpha$, consisting of a finite number of open arcs
 of the circle $|z|=1$,  into the circle $|w|=1$, 
then  the inequality (5.2) in view of (3.3) gives

\[
\sqrt{|f'(0)|}\,\mathrm{cap}\, \left\{w:|w|=1, \,\, w \notin
f(\alpha) \right\} \leq
\]
\[
\leq \left(\mathrm{cap}\,\{z{:}|z|{=}1, z \notin \alpha\} \right)
\exp \left\{-\sum\limits_{\nu\geq 1}\, \frac{n_\nu}{2}\, g_U(0,
z_{\nu}, \{z{:}|z|{=}1, z \notin \alpha\}) \right\},
\]
where $U=\{z:|z|<1\}$. 	In the case when $B=U,$  $\gamma=\partial U,$
$z_0=0$ from the inequality (5.2) there follows Hayman's inequality [H, p.
124]

\[
|f'(0)|\leq r(f(U), \, w_0).
\]
For an arbitrary domain
 $B$ and $\gamma=\partial B$ Corollary 5.2 was established
by Mityuk [M1]. We now give an example of the inequality (5.2),
when domains  $B$ and $G$ both are half planes. Let the function $f$ 
be regular in $\mathrm{Im}z>0$, $f(\infty)=\infty,$ $\mathrm{Im}
f(z)>0 $ for $\mathrm{Im} z>0$, and let $f$ map the set
$\{z=x+iy:|x|>1,  y=0\}$ into the set $\{w=u+iv:u\notin[a, b],
v=0\}.$ Then

\[
b-a\leq 2.
\]
In fact, from (5.2) for every point $\zeta$ of the upper half plane
we have

\[
r(\{z:\mathrm{Im}z>0\}, [-1, 1], \zeta)|f'(\zeta)|\leq
r(\{w:\mathrm{Im}w>0\}, [a, b], f(\zeta)).
\]
Letting $\zeta \to \infty$, we get

\[
r(\{z:\mathrm{Im}z>0\}, [-1, 1], \infty)\leq
r(\{w:\mathrm{Im}w>0\}, [a, b], \infty),
\]
which is equivalent to the required inequality.

\vspace{0.5cm} {\bf  5.3. Corollary.}  {\it Let the function} $f$
{\it be regular and }  $p$-{\it valent in the domain} $B:=\{z=x+iy:x>0,
y>0\},$ $f(B)\subset B,$ {\it and } $f(\alpha)\subset\alpha,$ {\it
where} $\alpha=\{z=x+iy:x>0, y=0\}.$ {\it We assume that in some
neighborhood of the point } $z_0 \in B$ 
{\it we have the expansion
}
\[
f(z)= z_0+c_p(z-z_0)^p+ \ldots
\]
  {\it
and let for some point} $\zeta \in B,$ {\it distinct from}
$z_0,$  {\it there exist exactly } $p$ 
{\it  distinct preimages in the domain} $B:z_j, j=1, \ldots, p.$ 
{\it  The the inequality}

\[
\left|\frac{\zeta\mathrm{Re}\zeta}{\mathrm{Im}\zeta} \right|^p\,\,
 \left[\zeta\atop z_0 \right]^{2p}\geq |c_p|
 \left\{\prod\limits_{j=1}^p\,\left|\frac{f'(z_j)z_j\mathrm{Re}z_j}{\mathrm{Im}z_j}
 \right|
  \left[z_j\atop z_0 \right]^{2}
 \right\}\prod\limits_{j, k=1\atop j\not=k}^p\,
 \left[z_j\atop z_k \right]^{-1} \quad (5.3)
\]
 {\it holds where the following notation is used}

 \[
\left[a\atop b \right]=
\left|\frac{(a-b)(a-\overline{b})}{(a+b)(a+\overline{b})} \right|.
\]
 {\it If the function } $f$ {\it
defines a complete} $p$--{\it valent covering of } $B$ {\it
onto the domain} $B$ {\it such that} $f(\alpha)=\alpha,$  $f((\partial
B\backslash \alpha))=(\partial B)\backslash\alpha,$ {\it then in the inequality}
 (5.3) {\it the equality sign holds}.

\vspace{0.5cm} {\bf  Proof.} By Theorem 4.3

\[
p\left\{\log r(B, (\partial B)\backslash\alpha, \zeta)-
2g_B(\zeta, z_0, (\partial B)\backslash\alpha) \right\} \geq
\sum\limits_{j=1}^p\, \log[|f'(z_j)|\, r(B,
\]

\[
(\partial B)\backslash\alpha,
z_j)]+\log|c_p|-2\sum\limits_{j=1}^p\, g_B(z_j, z_0,(\partial
B)\backslash\alpha)+\sum\limits_{j, k=1\atop j\not=k}^p\, g_B(z_j,
z_k,(\partial B)\backslash\alpha).
\]
This inequality coincides with (5.3), as we see by observing that for 
$D:=\{z:\mathrm{Re}z>0\}$ we have

\[
g_B(z, \zeta, (\partial B)\backslash\alpha)=g_D(z, \zeta)+g_D(z,
\overline{\zeta})=-\log\left[z\atop \zeta \right],
\]
and, therefore,

\[
\log r(B, (\partial B)\backslash\alpha,
\zeta)=\lim\limits_{z\to\zeta}[g_B(z, \zeta, (\partial
B)\backslash\alpha)+ \log|z-\zeta|]=
\]
\[
{\hspace{-9.0cm}=\log\left|\frac{2\zeta\mathrm{Re}\zeta}{\mathrm{Im}\zeta}\right|.}
\]


\vspace{1.5cm}\noindent {\Large\bf 6  Open problems}

\vspace{0.5cm} {\bf 6.1.} Prove two-point distortion theorems for functions
$f$ univalent in the disk $|z|<1,$ involving a lower estimate for the difference
$|f(z_1)-f(z_2)|$ in terms of the 
$|z_1-z_2|$ and the Maclaurin coefficients of  $f$.

\vspace{0.5cm} {\bf 6.2.} Known two-point distortion theorems in the disk with one point $z_1$
fixed and the other point 
$z_2$ approaching the boundary give in the limit either trivial estimates or
classical ones (provided that boundary derivative exists). 
Prove two-point distortion theorems which
give more interesting information when one of the points tends to the boundary.

\vspace{0.5cm} {\bf 6.3.} Prove a nontrivial two-point distortion theorem for functions
regular (not necessarily univalent) in the unit disk, involving the Schwarzian derivative.



 {\it }{\it }{\it }{\it }{\it
}{\it }{\it }{\it }{\it }{\it }{\it }{\it }{\it }{\it }{\it }{\it
}{\it }{\it }{\it }{\it }{\it }{\it }{\it }{\it }{\it }{\it }{\it
}{\it }{\it }{\it }{\it }{\it }{\it }{\it }{\it }{\it }{\it }{\it
}{\it }{\it }{\it }{\it }{\it }{\it }{\it }{\it }{\it }{\it }{\it
}{\it }{\it }{\it }{\it }{\it }{\it }{\it }{\it }{\it }{\it }{\it
}{\it }{\it }{\it }{\it }{\it }{\it }{\it}
\bigskip

\vspace{0.5cm}
\noindent{\large {\bf References}}

\vspace{0.5cm}

\small
\noindent[A]\quad \textsc{Yu. E. Alenitsyn}: 
{\it About univalent functions without
common values in multiplyconnected domains,} 
Trudy Mat. Inst. AN SSSR {\bf 94} (1968), 4-18.

\noindent[BS]\quad \textsc{S.~Bergman and M.~Schiffer}: 
{\it Kernel functions and elliptic differential equations 
in mathematical physics}, Academic Press, New York, 1953.

\noindent[DitSol] \quad \textsc{ B. Dittmar and A. Yu. Solynin}: 
Distortion of the hyperbolic Robin capacity 
under a conformal mapping, and extremal configurations. (Russian) 
Zap. Nauchn. Sem. S.-Peterburg. Otdel. Mat. Inst. Steklov. (POMI) 263 (2000), 
Anal. Teor. Chisel i Teor. Funkts. 16, 49--69, 238; translation in 
J. Math. Sci. (New York) 110 (2002), no. 6, 3058--3069.

\noindent[Dub1]\quad \textsc{V.~N.~Dubinin}: {\it Generalized
condensers and the asymptotics of their capacities under a degeneration of
some plates},(Russian) Zap.Nauchn.Sem.S.-Peterburg.Otdel.Mat.Inst.
Steklov. (POMI),  \textbf{302} (2003),
  38-51, translation in J.Math Sci. (N.Y.)
\textbf{129} (2005), no 3, 3835--3842. 

\noindent[Dub2]\quad \textsc{V.~N.~Dubinin}: {\it The majorization
principle for $p$-valent functions},
(Russian) Mat. Zametki \textbf{65}(1999), no.4,
533--541; translation in Math. Notes  \textbf{65}(1999), no. 3--4,
447--453.

\noindent[DubKos]\quad \textsc{V.~N.~Dubinin and E.V.~Kostyuchenko}:
{\it Extremal problems associated with $n$-fold symmetry 
in the theory of functions.} (Russian) Zap. Nauchn. 
Sem. S.-Peterburg. Otdel. Mat. Inst. 
Steklov. (POMI) 276 (2001), Anal. Teor. Chisel i Teor. 
Funkts. 17, 83--111, 350; 
translation in J. Math. Sci. (N. Y.) 118 (2003), no. 1, 4778--4794.

\noindent[DubPr]\quad \textsc{V.~N.~Dubinin and  E.G.~Prilepkina}:
{\it On  the preservation of}{\it generalized reduced modulus under some
geometric trans\-formations of domains in the plane}, Far
Eastern Math. J. \textbf{6}(2005), 1--2, 39--56.

\noindent[D1]\quad  \textsc{P.~Duren}: {\it Robin capacity},
Computational methods and Function Theory (CMFT'97)/ N. Papamichael,
St. Ruscheweyh and E.B. Saff (Eds.), World scientific Publishing Co.(1999) 177-190.

\noindent[D2]\quad  \textsc{P.~L.~Duren}: {\it Distortion in
certain conformal mappings of an annulus}, Michigan Math. J. 10(1963),
431--441.

\noindent[DPT]\quad \textsc{P.~Duren, J.~Pfaltzgraff, and E.~Thurman}:
{\it Physical interpretation and further properties of Robin capacity.}
 Algebra i Analiz 9 (1997), no. 3, 211--219; 
translation in St. Petersburg Math. J. 9 (1998), no. 3, 607--614.

\noindent[DS1]\quad \textsc{P.~Duren and M.~Schiffer}:
{\it Robin functions and energy functionals of multiply connected domains},
Pacific J.Math., \textbf{148}
(1991), 251-273.

\noindent[DS2]\quad \textsc{P.~Duren and M.~Schiffer}:
{\it Robin functions and distortion of capacity under conformal mapping},
Complex Variables Theory Appl.\textbf{21}(1993), 189-196.

\noindent[GH]\quad \textsc{D.~Gaier and W.K.~Hayman}: {\it On the
computation of modules of } {\it long quadrilaterals}, 
Constr. Approx., \textbf{7} (1991), 453-467.

\noindent[H]\quad \textsc{W.~K.~Hayman}:
{\it Multivalent functions}, Second ed., Cambridge 
Univ. Press, Cambridge, 1994.

\noindent[Her]\quad \textsc{J.~Hersch}:{\it On the reflection
principle and some elementary}{\it  ratios of conformal radii},
J. Anal. Math. \textbf{44} (1984/85), 251-268.

\noindent[Huck]\quad \textsc{F.~Huckemann}: 
{\it Extremal elements in certain classes of conformal 
mappings of an annulus}, Acta. Math. 118 (1967), 193--221.

\noindent[J]\quad \textsc{J.~A.~Jenkins}:
 {\it On two-point distortion theorems for bounded 
univalent regular functions},
Kodai. Math. J. 24 (2001), 3, 329-338.

\noindent[KMin]\quad \textsc{S.~Kim and D.~Minda}:
{\it Two-point distortion theorems for univalent functions}, 
Pacific J. Math., \textbf{163}(1994), 137-157.

\noindent[KrR]\quad \textsc{D.~Kraus and O.~Roth}:
{\it Weighted distortion in conformal mapping in enclidean, hyperbolic and
elliptic geometry}, Ann. Acad. Sci.Fenn. Math., \textbf{31}(2006), 111-130.

\noindent[MaMin]\quad \textsc{W.~Ma and D.~Minda}: {\it Two-point
distortion for univalent functions}, J. Comput. Appl. Math.,
\textbf{105} (1999), 385-392.

\noindent[Mik]\quad \textsc{V.M. Miklyukov}: 
{\it Conformal mapping of irregular surface
and its applications.} Izdat. VolGU, Volgograd, 2005.

\noindent[M1]\quad \textsc{I.P.~Mityuk}: {\it Symmetrization principle
for multiply connected domains}, (Russian), Dokl. Akad. Nauk SSSR 
\textbf{157} (1964), 268-270.

\noindent[M2]\quad \textsc{I.P.~Mityuk}: {\it Symmetrization principle
for multiply connected domains and some of its applications}, (Russian),
Ukr. Mat. Zh. \textbf{17} (1965),4, 46-54.

\noindent[Nas]\quad \textsc{S.~Nasyrov}:
{\it Robin capacity and lift of infinitely thin airfoils},
 Complex Variables Theory Appl.\textbf{47}(2002), 2, 93-107.

\noindent[Neh]\quad \textsc{Z.~Nehari}: {\it Some inequalities in
the theory of functions},
 Trans. of Amer. Math. Soc., \textbf{75}(1953), no.2, 256-286.

\noindent[NT]\quad \textsc{M.~D.~O'Neill and R.~E.~Thurman}:
{\it Extremal domains for Robin capacity}, Complex Variables
\textbf{41}(2000),  91-109.

\noindent[P]\quad \textsc{Ch.~Pommerenke}: {\it Boundary behaviour
of conformal maps}, Berlin: Springer, 1992.

\noindent[ProW]\quad \textsc{M.~H.~Protter and H.F. Weinberger}: {\it
Maximum principles in differential equations}, Springer-Verlag,
New York,Inc., 1984.

\noindent[S]\quad \textsc{M.~Schiffer}: {\it
Some recent developments in the theory of conformal mapping}, 
in Dirichlet's Principle,Conformal
Mapping and Minimal Surfaces (R. Courant: Interscience:
New York, 1950), appendix.

\noindent[Sin]\quad \textsc{V.~Singh}:
{\it Grunsky inequalities and coefficients of bounded schlicht functions},
 Ann. Acad. Sci. Fenn.  \textbf{310}(1962), 1-21.

\noindent[Sol1]\quad \textsc{
A. Yu.  Solynin}:  Moduli and extremal metric problems. 
(Russian) Algebra i Analiz 11 (1999), no. 1, 3--86; 
translation in St. Petersburg Math. J. 11 (2000), no. 1, 1--65.

\noindent[Sol2] \quad \textsc{A. Yu. Solynin}: 
Boundary distortion and extremal problems in some 
classes of univalent functions. (Russian) Zap. Nauchn. Sem. 
S.-Peterburg. Otdel. Mat. Inst. Steklov. (POMI) 204 (1993), 
Anal. Teor. Chisel i Teor. Funktsii. 11, 115--142, 169--170; 
translation in J. Math. Sci. 79 (1996), no. 5, 1341--1358. 

\noindent[St1]\quad \textsc{M.~Stiemer}:
{\it A representation formula for the Robin function},
 Complex Variables Theory Appl.\textbf{48}(2003), 5, 417-427.

 \noindent[St2]\quad \textsc{M.~Stiemer}:
{\it Extremal point methods  for  Robin capacity},
 Comp. Meth. Funct. Th. \textbf{4}(2004), 2, 475-496.

 \noindent[Sto]\quad \textsc{S.~Stoilow}:
{\it Teoria functiilor de o variabila complexa}, vol. II,
 Editura Acad. Rep. Pop. Romine, 1958.

\noindent[Vas]\quad \textsc{A.~Yu.~Vasil'ev}:
{\it  Robin's modulus in a Hele-Shaw problem},
 Complex Variables Theory Appl.\textbf{49}(2004), 7--9, 663-672.

\noindent[V]\quad \textsc{
 M. Vuorinen}: {\em Conformal Geometry
and Quasiregular Mappings}, Lecture Notes in Math. 1319,
Springer-Verlag, Berlin--New York, 1988.

\bigskip
 \noindent
{\bf V.N. Dubinin }\\
Institute of Applied Mathematics\\
Far East Branch Russian Academy of Sciences\\
Vladivostok\\
RUSSIA\\
Email: {\tt dubinin@iam.dvo.ru}\\
 \medskip

\noindent
{\bf M. Vuorinen}\\
Department of Mathematics\\
FIN-20014 University of Turku \\
FINLAND\\
E-mail: {\tt vuorinen@utu.fi}\\
\end{document}